\documentclass{article}
\usepackage[utf8]{inputenc}

\usepackage{amsmath}
\usepackage{amssymb}
\usepackage{mathtools}
\usepackage{graphicx}

\title{There is a $3\times3$ Magic Square of Squares on the Moon -- A Lot of Them, Actually}
\author{Christian Woll \\
cwoll@ucsd.edu}
\date{November 2018}

\begin{document}

\maketitle

\section{Lunar Arithmetic}
Magic squares have been well explored here on Earth$^{[1]}$, but there appears to have been little-to-no examination of magic squares on the Moon.\newline
\indent Lunar arithmetic (originally called ``dismal" arithmetic) is an alternative system of adding and multiplying numbers$^{[2]}$. To briefly cover the definitions given in the original paper, adding two numbers is done by taking the larger of each digit, 
$$15+83=85.$$ Multiplication is done by taking the smaller of each digit,
$$17\times 3=13, \quad 5\times 5=5.$$
One performs traditional distribution where both numbers are multiple digits long. Before continuing, readers not yet familiar with lunar arithmetic are strongly encouraged to skim the original paper$^{[2]}$ and/or watch a recent video in which Neil Sloane himself explains the number system$^{[6]}$.\\
\indent Of course, lunar arithmetic can be done in any base, not just base $10$. One may think of base $B$ lunar arithmetic as a the semi-ring of polynomials $\{0,1,...,B-1\}[x]$ where coefficient addition and multiplication are the $\textit{max}$ and $\textit{min}$ functions respectively. \\
\indent We will explore here the properties of $3\times3$ magic squares under lunar arithmetic. Remember that a magic square is a grid of integers in which the entries of each row, column, and the two diagonals sum to the same total.
\begin{equation*}
\text{Take}\quad
\begin{tabular}{c|c|c}
  $8$ & $1$ & $6$ \\      \hline
  $3$ & $5$ & $7$ \\      \hline
  $4$ & $9$ & $2$
\end{tabular}\quad \text{for example.}
\end{equation*}
Sometimes one restricts the entries of the square to be the integers $1$ through $n^2$ (where $n=3$ in the former case). Here we disregard this restriction but will still attempt to find squares with distinct entries. We don't care about trivial solutions like 
\begin{equation*}
\begin{tabular}{c|c|c}
  $42$ & $42$ & $42$ \\      \hline
  $42$ & $42$ & $42$ \\      \hline
  $42$ & $42$ & $42$
\end{tabular}\ .
\end{equation*}

\section{General Magic Squares}
There are, in fact, magic squares on the Moon (again, throughout the paper we say ``on the Moon" to mean ``under dismal/lunar arithmetic" and unless stated otherwise, it is assumed from here on that any mention of a magic square refers to a $\textit{lunar}$ magic square). Consider the example 
\begin{equation*}
\begin{tabular}{c|c|c}
  $12$ & $0$ & $20$ \\      \hline
  $1$ & $22$ & $10$ \\      \hline
  $21$ & $20$ & $2$
\end{tabular}\ .
\end{equation*}
Is it magic? Let's check the sums:
$$12+0+20\ =\ 1+22+10\ =\ 21+20+2\ =\ 22$$
$$12+1+21\ =\ 0+22+20\ =\ 20+10+2\ =\ 22$$
$$12+22+2\ =\ 21+22+20\ =\ 22$$
It's magic -- the sums are all $22$. It is claimed this is also the smallest possible total of a magic square in any base $B> 2$. The proof is left undone.\\
\indent But does there exist a binary ($B=2$) lunar magic square? Yes:
\begin{equation*}
\begin{tabular}{c|c|c}
  $1111$ & $1110$ & $1011$ \\      \hline
  $1010$ & $0$ & $111$ \\      \hline
  $110$ & $1001$ & $1$
\end{tabular}\ 
\end{equation*}
This is proof by construction that lunar magic squares exist in all bases, $B\ge 2$. \\
\indent Before going on, it is best to introduce the notion of one number $\textit{dominating}$ another. We say a lunar integer, $m$, dominates another, $n$, if the digits of $m$ are greater than or equal to the digits of $n$ paired respectively. This is equivalent to $m+n=m$. For example, $287$ dominates $185$. We use the notation
$$m \gg_B n, \quad 287 \gg_B 185$$
to denote domination in base $B$. \\
\indent In both of the lunar squares we have seen so far, there is a single entry which dominates all of the others ($22$ and $1111$ respectively). Note that for an entry to dominate all others it must be the total. \\
\indent Are there lunar magic squares in which no individual entry is also the total? Yes: 
\begin{equation*}
\begin{tabular}{c|c|c}
  $40$ & $34$ & $41$ \\      \hline
  $42$ & $0$ & $24$ \\      \hline
  $14$ & $43$ & $4$
\end{tabular}\ .
\end{equation*}
The total is $44$ but the largest entry is $43$. This square brings us to another important observation. Since there are no carries in lunar arithmetic, a magic square is magic for any subset of the digits. For example, if we take only the $2$nd digit (counting from the right) of each entry in the former square, we get
\begin{equation*}
\begin{tabular}{c|c|c}
  $4$ & $3$ & $4$ \\      \hline
  $4$ & $0$ & $2$ \\      \hline
  $1$ & $4$ & $0$
\end{tabular}\ .
\end{equation*}
The square no longer has distinct entries -- but what concerns us is that its sums are still all the correct. And to be clear, we are imagining each entry as buffered on the left by zeros, $3=...0003$. Thus we say that the $2$nd digit of $3$ is $0$ (just as the $2$nd digit of $2745$ is $4$).\\
\indent It follows that any magic square can be represented as an element-wise sum of single digit magic squares,
\begin{equation}
\begin{tabular}{c|c|c}
\label{base10}
  $40$ & $34$ & $41$ \\      \hline
  $42$ & $0$ & $24$ \\      \hline
  $14$ & $43$ & $4$
\end{tabular}\ =\ 
\begin{tabular}{c|c|c}
  $40$ & $30$ & $40$ \\      \hline
  $40$ & $0$ & $20$ \\      \hline
  $10$ & $40$ & $0$
\end{tabular}\ +\ 
\begin{tabular}{c|c|c}
  $0$ & $4$ & $1$ \\      \hline
  $2$ & $0$ & $4$ \\      \hline
  $4$ & $3$ & $4$
\end{tabular}\ .
\end{equation}
\indent So we ask: in any single-digit magic square, what are the minimal entries which must be included to cover all $8$ sums? It turns out there are only two arrangements (not counting rotations or reflections as distinct):
\begin{equation*}
\begin{tabular}{c|c|c}
  $1$ & $0$ & $1$ \\      \hline
  $1$ & $0$ & $0$ \\      \hline
  $0$ & $1$ & $0$
\end{tabular}\quad \text{and}\quad 
\begin{tabular}{c|c|c}
  $1$ & $0$ & $0$ \\      \hline
  $0$ & $1$ & $0$ \\      \hline
  $0$ & $0$ & $1$
\end{tabular}\ .
\end{equation*}
There are in fact other single-digit arrangements that cover all the sums, like
\begin{equation*}
\begin{tabular}{c|c|c}
  $1$ & $1$ & $1$ \\      \hline
  $0$ & $0$ & $1$ \\      \hline
  $0$ & $1$ & $1$
\end{tabular}\ ,
\end{equation*}
but by ``minimal entries" we mean that if any entry is removed (i.e. set to $0$), the square ceases to be magic. \\
\indent Manufacturing magic squares is easy given this information. For example, if we choose any $a\gg_B b,c,d$ and $\alpha \gg_B \beta, \gamma, \delta$, then a square may be produced as follows:
\begin{equation*}
\begin{tabular}{c|c|c}
  $a0...0$ & $b0...0$ & $a0...0$ \\      \hline
  $a0...0$ & $0$ & $c0...0$ \\      \hline
  $d0...0$ & $a0...0$ & $0$
\end{tabular}
\ +\ 
\begin{tabular}{c|c|c}
  $0$ & $\alpha$ & $\delta$ \\      \hline
  $\gamma$ & $0$ & $\alpha$ \\      \hline
  $\alpha$ & $\beta$ & $\alpha$
\end{tabular}
\ =\ 
\begin{tabular}{c|c|c}
  $a0...0$ & $b\alpha$ & $a\delta$ \\      \hline
  $a\gamma$ & $0$ & $c\alpha$ \\      \hline
  $d\alpha$ & $a\beta$ & $\alpha$
\end{tabular}\ .
\end{equation*}
The total is $a\alpha$. For example, the square in equation \ref{base10} is generated by choosing $(a,b,c,d)=(\alpha,\beta,\gamma,\delta)=(4,3,2,1)$. 

\section{Magic Squares of Powers}
It is unknown whether or not a $3\times3$ magic square of squared integers exists on Earth$^{[5]}$. It turns out that an abundance of magic squares of squares exist on the Moon. \\
\indent We should define what a lunar square is before going on. It is simply any lunar integer multiplied by itself ($1^2=1,2^2=2,...,57^2=557,...$ $^{[4]}$). Our first example can be obtained with the parametrization of in the previous section:
\begin{equation*}
\begin{tabular}{c|c|c}
  $440$ & $330$ & $440$ \\      \hline
  $440$ & $0$ & $220$ \\      \hline
  $110$ & $440$ & $0$
\end{tabular}
\ +\ 
\begin{tabular}{c|c|c}
  $4$ & $8$ & $5$ \\      \hline
  $6$ & $0$ & $8$ \\      \hline
  $8$ & $7$ & $8$
\end{tabular}
\ =\ 
\begin{tabular}{c|c|c}
  $444$ & $338$ & $445$ \\      \hline
  $446$ & $0$ & $228$ \\      \hline
  $118$ & $447$ & $8$
\end{tabular}\ 
\end{equation*}
\begin{equation*}
\ =\ 
\begin{tabular}{c|c|c}
  $44^2$ & $38^2$ & $45^2$ \\      \hline
  $46^2$ & $0^2$ & $28^2$ \\      \hline
  $18^2$ & $47^2$ & $8^2$
\end{tabular}\ .
\end{equation*}
The total is $48^2=448$. But there is a magic square of squares with a smaller total, $24^2$:
\begin{equation*}
\begin{tabular}{c|c|c}
  $22^2$ & $0^2$ & $14^2$ \\      \hline
  $1^2$ & $24^2$ & $2^2$ \\      \hline
  $4^2$ & $3^2$ & $23^2$
\end{tabular}\ .
\end{equation*}
This seems to be the smallest total of any magic square of squares.

\indent It was remarked without proof in the original paper that lunar integers whose digits are non-decreasing (eg. $1134448$) are closed under addition and multiplication$^{[2]}$ (see the end of section 2). We give here, without proof the formula for the $n$th power of such an integer. If $a=a_k...a_1a_0|_B$ and $a_{i+1} \le a_i$ then
$$a^n=\overbrace{a_k...a_k}^{\text{$n$ times}}...\overbrace{a_1...a_1}^{\text{$n$ times}}a_0.$$
For example, $1134448^3=1111113334444444448$. It follows immediately that our magic square of squares is valid for any alternative power. That is to say, 
\begin{equation*}
\begin{tabular}{c|c|c}
  $44^3$ & $38^3$ & $45^3$ \\      \hline
  $46^3$ & $0^3$ & $28^3$ \\      \hline
  $18^3$ & $47^3$ & $8^3$
\end{tabular}\ ,\quad 
\begin{tabular}{c|c|c}
  $44^4$ & $38^4$ & $45^4$ \\      \hline
  $46^4$ & $0^4$ & $28^4$ \\      \hline
  $18^4$ & $47^4$ & $8^4$
\end{tabular}\ ,\quad 
\begin{tabular}{c|c|c}
  $44^5$ & $38^5$ & $45^5$ \\      \hline
  $46^5$ & $0^5$ & $28^5$ \\      \hline
  $18^5$ & $47^5$ & $8^5$
\end{tabular}\ ,\quad...
\end{equation*}
are also magic squares. This may be interpreted as proof that infinitely many lunar magic squares of any power exist. There are many other such families as well. Consider 
\begin{equation*}
\begin{tabular}{c|c|c}
  $1447^n$ & $1347^n$ & $1444^n$ \\      \hline
  $1446^n$ & $0^n$ & $1247^n$ \\      \hline
  $1147^n$ & $1445^n$ & $1^n$
\end{tabular}\ .
\end{equation*}
\indent But, is there a lunar magic square of squares in base $2$? Yes: 
\begin{equation}
\label{base2}
\begin{tabular}{c|c|c}
  $11^2$ & $101^2$ & $1001^2$ \\      \hline
  $110^2$ & $1011^2$ & $1^2$ \\      \hline
  $1010^2$ & $0^2$ & $111^2$
\end{tabular}\ =\ 
\begin{tabular}{c|c|c}
  $111$ & $10101$ & $1001001$ \\      \hline
  $11100$ & $1011111$ & $1$ \\      \hline
  $1010100$ & $0$ & $11111$
\end{tabular}\ .
\end{equation}
\indent Each entry in this binary square has $7$ digits. This is also the smallest possible total of any binary magic square of squares.\\
\indent In all the examples of magic squares we have seen, the total has been a lunar square ($48^2$, $24^2$, and $1011^2$). So does there exist a magic square of squares whose total is not a lunar square? Yes:
\begin{equation*}
\begin{tabular}{c|c|c}
  $39^2$ & $40^2$ & $29^2$ \\      \hline
  $19^2$ & $33^2$ & $41^2$ \\      \hline
  $42^2$ & $9^2$ & $43^2$
\end{tabular}\ .
\end{equation*}
The total, $439$, is not a lunar square (which can be checked by brute force since there are only $100$ lunar squares with $3$ digits or less). In fact, we have found a magic square of squares whose total is a lunar prime(!). See section 3 of the original paper to find out what exactly that means$^{[2]}$.

\section{Pythagorean Triples}
There are, in fact, Pythagorean triples on the Moon:
$$22^2+4^2=24^2.$$
Many of these triples are trivial (example: $3^2+3^2=3^2$). So to keep things interesting, we will define a lunar Pythagorean triples to be \textit{distinct}.\\
\indent It follows from the observations we made in Section 2 that the element-wise sum of any two magic squares is also a magic square (just like back on Earth!). So do there exist magic squares of squares whose element-wise sum is also a magic square of squares? Yes:
\begin{equation*}
\begin{tabular}{c|c|c}
  $5789^2$ & $5778^2$ & $6778^2$ \\      \hline
  $6678^2$ & $6789^2$ & $5689^2$ \\      \hline
  $6788^2$ & $5678^2$ & $6689^2$
\end{tabular}\ +\ 
\begin{tabular}{c|c|c}
  $13458^2$ & $13348^2$ & $23348^2$ \\      \hline
  $22348^2$ & $23458^2$ & $12458^2$ \\      \hline
  $23448^2$ & $12348^2$ & $22458^2$
\end{tabular}
\end{equation*}
\begin{equation*}
=\ 
\begin{tabular}{c|c|c}
  $15789^2$ & $15778^2$ & $26778^2$ \\      \hline
  $26678^2$ & $26789^2$ & $15689^2$ \\      \hline
  $26788^2$ & $15678^2$ & $26689^2$
\end{tabular}\ .
\end{equation*}
There are nine Pythagorean triples represented here.

\section{Questions}
1) What is the smallest $3\times3\times3$ magic $\textit{cube}$ of lunar squares in base $2$?  (Note that ``smallest" is open to interpretation as the lunar integers have no ``fully satisfactory" ordering$^{[2]}$; see Section 5 of the original paper). The smallest such \textit{square} was given in equation \ref{base2}.\\
2) What is the smallest example of two magic squares of squares whose element-wise sum is also a magic square of squares in base $2$?
\pagebreak

$$\textbf{REFERENCES}$$\\\
[1] http://www.multimagie.com/English/Enigmas.htm\\\\\
[2] Applegate, LeBrun, \& Sloane, $\textit{Dismal Arithmetic}$\\\\\
[3] SageMath, the Sage Mathematics Software System (Version 6.9),\\
\indent The Sage Developers, 2015, http://www.sagemath.org.\\\\\
[4] http://oeis.org/A087019\\\\\
[5] http://www.multimagie.com/English/SquaresOfSquaresSearch.htm\\\\\
[6] Sloane, \textit{Primes on the Moon (Lunar Arithmetic)}, Numberphile,\\
\indent $[$interview by Brady Haran$],$ \\
\indent https://www.youtube.com/watch?v=cZkGeR9CWbk\\\\\

\end{document}